\documentclass[12pt]{article}
\usepackage{amsmath}
\usepackage{amscd}
\usepackage{amssymb}
\usepackage{array}

\newtheorem{definition}{Definition}[section]
\newtheorem{theorem}{Theorem}[section]
\newtheorem{example}{Example}[section]
\newtheorem{proposition}{Proposition}[section]
\newtheorem{lemma}{Lemma}[section]

\begin{document}

\newcommand{\R}{\mathbb{R}}
\newcommand{\C}{\mathbb{C}}
\newcommand{\g}{\mathcal{G}}
\newcommand{\A}{\mathcal{A}}
\newcommand{\I}{\mathcal{I}}
\newcommand{\m}{\mathbf{m}}
\newcommand{\La}{\mathcal{L}}
\newcommand{\N}{\mathcal{N}_{\Theta^r}}
\newcommand{\de}{\delta}
\newcommand{\Z}{\mathbb{Z}}
\newcommand{\dex}{\partial_x}
\newcommand{\dey}{\partial_y}
\vskip 2cm
  \centerline{\LARGE \bf Algebraic and Differential  Star Products   }
\bigskip
\centerline {\LARGE\bf on  Regular Orbits of Compact Lie Groups}

\vskip 2.5cm

\centerline{R. Fioresi$^\star$\footnote{Investigation supported by
the University of Bologna, funds for selected research topics.},
A. Levrero$^\ast$ and M. A. Lled\'o$^\dagger$}

\bigskip

\centerline{\it $^\star$Dipartimento di Matematica, Universit\`a
di Bologna }
 \centerline{\it Piazza di Porta S. Donato, 5.}
 \centerline{\it40127 Bologna. Italy.}
\centerline{{\footnotesize e-mail: fioresi@dm.UniBo.it}}

\bigskip

\centerline{\it $^\ast$ Dipartimento di Fisica, Universit\`a di
Genova. } \centerline{\it Via Dodecaneso 33, 16146 Genova, Italy,
and  } \centerline{\it INFN, Sezione di Genova. Italy.}
\centerline{{\footnotesize e-mail: levrero@ge.infn.it}}

\bigskip

\centerline{\it $^\dagger$ Dipartimento di Fisica, Politecnico di
Torino,} \centerline{\it Corso Duca degli Abruzzi 24, I-10129
Torino, Italy, and} \centerline{\it INFN, Sezione di Torino,
Italy.} \centerline{{\footnotesize e-mail:
lledo@athena.polito.it}}
 \vskip 2cm

\begin{abstract}

In this paper we study  a family of algebraic deformations of
regular coadjoint orbits of compact semisimple Lie groups with the
Kirillov Poisson bracket. The deformations are restrictions of
deformations on the dual of the Lie algebra.  We prove that there
are non isomorphic deformations in the family. The star products
are not differential, unlike the star products considered in other
approaches. We make a comparison with the differential star
product canonically defined by Kontsevich's map.
\end{abstract}

\vfill\eject

\section{Introduction}

Coadjoint orbits of Lie groups are symplectic manifolds that can
be used to model physical systems that have a continuous group of
symmetries. The Kirillov-Kostant orbit principle allows in many
cases to associate canonically a unitary representation to the
orbit. The Hilbert space of the representation can then be thought
as the Hilbert space of the quantum theory. A quantization map
which takes a class  of functions on the phase space to operators
in such Hilbert space can be constructed.  This is the approach of
geometric quantization (see Ref. \cite{vo} for a review).

 On the other hand, the pioneering work by Bayen {\it et al.} \cite{bffls} on
deformation quantization  raised the problem of quantizing
 the coadjoint orbits with a
 radically different method. However,  being based on the same physical principles,  it is
natural to expect a relation between the two approaches. In fact,
it was thought that deformation quantization, which ``forgets"
about the Hilbert space on which the quantum algebra is
represented, could nevertheless throw light on the
Kirillov-Kostant orbit principle \cite{fr}. The  algebra that
appears in geometric quantization is defined as the quotient of
the enveloping algebra by a prime ideal which is contained in the
kernel of the corresponding representation  \cite{vo}. The method
of geometric quantization is however more general than the
Kirillov-Kostant orbit principle. A comparison  with deformation
quantization for the case of $\R^{2n}$ with the standard
symplectic structure was done in Ref. \cite{gv}.

In the work of Bayen  {\it et al.} \cite{bffls} only flat
symplectic manifolds were studied. The existence of a deformation
quantization of general symplectic manifolds was first established
by De Wilde and LeComte \cite{dl}, and using different methods by
Omori, Maeda and Yoshioka \cite{omy} and by Fedosov \cite{fe}. For
a comparison between the methods of De Wilde and Le Comte and
Fedosov, see Ref. \cite{de}. In Ref. \cite{ma}, the existence of
tangential deformations for any regular Poisson manifold was
proven. In Ref. \cite{ko}, Kontsevich settled the fundamental
question of the existence of deformations for arbitrary (formal)
Poisson manifolds. In all these works the deformations are taken
to be differential, that is, the product structure
  in the deformed algebra is defined through bidifferential
  operators.

Explicit star products for non flat manifolds are not easy to
construct. In Ref. \cite{gu}, Gutt constructed a star product on
the cotangent bundle of a Lie group. In Ref. \cite{cg}, Cahen and
Gutt constructed a deformation of the algebra of polynomials on
the regular coadjoint orbits of compact semisimple groups, using
the fact that the universal enveloping algebra is  a deformation
of the algebra of polynomials on the dual of the Lie algebra
\cite{ho}. They showed that, although the deformation  on the
whole space is differential, the one induced on the orbit is not.
Moreover, in Ref. \cite{cgr} they show that for semisimple groups
``tangential" deformations (that is, deformations
  on the ambient space that  restrict well to the orbits)
that are at the same time differential and that extend over the
origin do not exist. Deformations of coadjoint orbits  were also
studied in \cite{alm} in terms of a polarization of the orbit
(also used in geometric quantization). The resulting star product
is covariant. More generally, deformations of K\"ahler manifolds
were studied in \cite{cgr2}

In Ref. \cite{fl} a family of star products on coadjoint orbits of
semisimple Lie groups was constructed as a quotient of   the
enveloping algebra by a suitable ideal. With a certain choice
inside the family of deformations one
  obtains the same star product as in Ref \cite{cg}. For another choice,
  in the special case of SU(2), the deformed algebra turns out to
  be the one of geometric quantization \cite{fl}.  In this case we
  can associate to the deformation quantization a unitary  representation
  in the spirit of  Berezin \cite{be}.

In the present work we further study the properties of this family
of deformations.  The organization of the paper is as follows. In
Section 2. we review the construction of the algebraic star
products on the orbit \cite{fl} and show that there exist non
equivalent products associated with a given algebraic Poisson
bracket.  We  also show that the ideal used to
 quotient the
enveloping algebra is prime. In fact, in geometric quantization
 the quantum algebra is the enveloping algebra
modulo a prime ideal; this ideal  is contained in the kernel of
the representation. In Section 3. we use Kontsevich's theorem on
differential star products to show that the  Kirillov Poisson
structure on the dual space of the Lie algebra of a semisimple Lie
group has only one possible deformation.  In Section 4. we study
the algebraic star products on the orbit and show that they are
not, in general, differential. In Section 5. we show different
ways of constructing star products on the orbit.

\section{Deformation of the polynomial algebra of a regular orbit}

In this section we review the results of Ref.  \cite{fl} where a
family of  deformations of the polynomial algebra of a regular
coadjoint orbit of a semisimple Lie group was constructed. We show
that the different deformations in the family are not necessarily
equivalent by exhibiting a counterexample.

Let $G$ be  a  complex semisimple Lie group of dimension $n$ and
rank $m$,    $\g$ its Lie algebra and $U$ the enveloping algebra
of $\g$. Let $T_{\C}(\g)$ be the full tensor algebra of $\g$ over
$\C$. Consider the  algebra $T_{\C}(\g)[[h]]$ and its proper two
sided ideal \begin{equation} \mathcal{L}_{[h]}=\sum_{X,Y \in \g}
T_{\C}(\g)[[h]] \otimes(X \otimes Y - Y \otimes X - h[X,Y])
\otimes T_{\C}(\g)[[h]]. \label{lh}\end{equation} We define
$U_{[h]}= T_{\C}(\g)[[h]]/\mathcal{L}_{[h]}$.

\begin{definition}
An associative algebra $\A_{[h]}$ over $\C[[h]]$ is a formal deformation of
a Poisson  algebra $(\A, \{\;,\;\})$ over $\C$ if there exists an
isomorphism of $\C[[h]]$-modules $\psi: \A[[h]]\longrightarrow
\A_{[h]} $ satisfying the following properties

\noindent {\bf a.} $\psi(f_1f_2)=\psi(f_1)\psi(f_2)$ mod($h$).

\noindent {\bf b.}
$\psi(f_1)\psi(f_2)-\psi(f_2)\psi(f_1)=h\psi(\{f_1,f_2\})$
mod($h^2$). \label{fd}
\end{definition}

Because of its relation with the problem of quantization, $\A_{[h]}$
is sometimes called a {\it deformation quantization} of $\A$.

 Notice that in the above definition we can substitute $\C[[h]]$ by $\C[h]$.
The  algebra is then a  module over $\C[h]$ which will be denoted
as $\A_h$. We will say that $\A_h$ is a {\it $\C[h]$-deformation}
of $\A$. Notice that a $\C[h]$-deformation extends to a formal
deformation, but the converse is not always true. Also, a
$\C[h]$-deformation can be specialized to any value of the
parameter $h$, since the ideal generated by the element $h-h_0$ is
proper in $\A_h$. One obtains then a complex algebra, $\A_{h_0}=
\A_h/(h-h_0)$.

It is well known that $U_{[h]}$ is a formal deformation of
$\C[\g^*][[h]]$ equipped with the Kirillov Poisson bracket
\cite{ho}.

We denote by $p_i,\; i=1,\dots m$ the algebraically independent
homogeneous  generators of the subalgebra of invariant polynomials
on $\g^*$,
\begin{equation}
I=\{p\in \C[\g^*] \; | \;  p(\hbox{Ad}^*(g)\xi)=p(\xi) \quad
\forall \xi \in \g^*, \, g \in G\}=\C[p_1,\dots p_m],
\label{invpol}
\end{equation}
 given by Chevalley's theorem.
If $S(\g)$ is the algebra of symmetric tensors on $\g$, we can
identify canonically Pol$(\g^*)=\C[\g^*]\approx S(\g)$.

Let $\{X_1,\dots , X_n\}$ be a basis for $\g$ and let $\{x_1,\dots
, x_n\}$ be the corresponding
generators of  $\C[\g^*]$. Then the symmetrizer
map $\hbox{Sym} :\C[\g^*]\longrightarrow  T_{\C}(\g)$ is given by
\begin{equation} \hbox{Sym}(x_1\cdots x_p)=\frac{1}{p!}\sum_{s \in
S_p}X_{s(1)}\otimes\cdots\otimes X_{s(p)} \label{sym}
\end{equation} where $S_p$ is the group of permutations of order
$p$. The composition of the symmetrizer with the natural
projection $T_{\C}(\g)\longrightarrow U$ is a linear isomorphism
that gives the identification $\C[\g^*]\approx U$. Moreover, it
sends the invariant polynomials $I$ isomorphically into the center
of $U$ (see for example Ref. \cite{va}). We can extend the
symmetrizer map as Sym: $\C[\g^*][[h]]\longrightarrow
T_{\C}(\g)[[h]]$, and   the projection
$\pi_h:T_{\C}(\g)[[h]]\longrightarrow U_{[h]}$. Then $P_i=\pi_h
\circ\mbox{Sym}(p_i)$ are also central elements. Note that $\pi_h
\circ\mbox{Sym}$ can be used as the isomorphism $\psi$ in
Definition \ref{fd}, $\psi:\C[\g]\rightarrow U_{[h]}$.

\bigskip

We  consider now  the  compact real form of $G$ with $\g^r$ the
real Lie algebra ($\g$ still denotes the complex Lie algebra). The
coadjoint orbits are algebraic manifolds given by
 the constrains,
$$ p_i(x_1,\dots x_n)=c^0_i,\quad c^0_i\in \R,\quad i=1,\dots m.
$$ There is a one to one correspondence from  the set of orbits
with the elements of a Weyl chamber in the Cartan subalgebra.
 Regular orbits are the orbits of elements in
the interior of a Weyl chamber, and they have maximal dimension.
Non regular orbits are given by constants $c^0_i$ satisfying some
constrains. This means that if $c^0_i$ define a regular orbit,
there is a neighborhood of the orbit that is  foliated with
regular orbits. We will use this property in the next  section.

Let $\Theta^r$ be a regular orbit. Then, the ideal of polynomials
in $\R[{\g^r}^*]$ that vanish on $\Theta^r$ is generated by the
elements $p_i-c_i^0$ \cite{kos}, so we can define $$ \I_0=
(p_i-c^0_i,\; i=1,...,m) \subset \R[{\g^r}^*], $$ and the algebra
of restrictions of polynomials to the orbit is
$\R[\Theta^r]=\R[{\g^r}^*]/\I_0$. We take the complexification of
this algebra $\C[\g^*]/\I_0$ (we denote still by $\I_0$ the ideal
in the complexified algebra), which is  the algebra of polynomials
on the complex orbit $\Theta$, $\C[\Theta]$. Consider a regular
orbit and define the two sided ideal in $U_{[h]}$ generated by the
elements $P_i-c_i(h)\; i=1,...,m$
\begin{equation} \I_{[h]}= (P_i-c_i(h), \;
i=1,...,m) \subset U_{[h]},\label{prime}
\end{equation}
where $c_i(h)\in\C[[h]]$ is such that $c_i(0)=c_i^0$ and $P_i={\rm
Sym}(p_i)$. In \cite{fl} it was shown the following
\begin{theorem}
The algebra $U_{[h]}/\I_{[h]}$  is a formal deformation of
$\C[\Theta]=\C[\g^*]/\I_0$.  $U_{h}/\I_{h}$  is a
$\C[h]$-deformation of $\C[\Theta]=\C[\g^*]/\I_0$.\label{flt}
\end{theorem}
Regularity is a technical assumption to show that
$U_{[h]}/\I_{[h]}$ is a free module isomorphic to
 $\C[\Theta][[h]].$

 \bigskip

The ideal $\I_0\subset \C[\g^*]$ is prime since the corresponding
algebraic variety is irreducible. We want to show now that the
ideal $\I_h\in U_h$ is prime. We define first a grading in $U_h$.
If $\{X_i, i=1,\dots, n\}$ is  a basis of $\g$ we set
$\deg(X_i)=\deg(h)=1$. This is a set of generators for
$T_{\C}(\g)[h]$. Notice that the relations in $\La_h$ (\ref{lh})
are homogeneous with respect to this grading, so a grading is
defined on $U_h$. The degree of an inhomogeneous element in $U_h$
is the maximal degree occurring in all of its monomials. Let us
restrict to modules over $\C[h]$.

\begin{proposition} Assume that $\deg(c_i(h))\leq \deg(P_i)$. Then if $FG\in \I_h$,
either $F\in \I_h$ or $G\in \I_h$. Hence $\I_h$ is prime.
\label{primet}\end{proposition}

\noindent {\it Proof.} Consider first the projection
$\rho:U_h\rightarrow U_h/hU_h\approx \C[\g^*]$. One has that for
any $F\in U_h$, $\deg(\rho(F))\leq\deg(F)$.

Since $\rho(P_i)=p_i$, $i=1,...,m$, we have $\rho(\I_h)=\I_0$.
Since $\deg(c_i(h))\leq \deg(P_i)$, if $f\in\I_0$ there exists
$F'\in\I_h$ with $\rho(F')=f$ and $\deg(F')=\deg(f)$. If
$f=f^i(p_i-c^i_0)$ one can take for example $F=\pi_h\circ{\rm
Sym}(f^i)(P_i-c^i(h))$.

 Assume that $FG\in \I_h$. Then $$\rho(FG)=\rho(F)\rho(G)=:fg\in \I_0,$$
where we denote the projections by small case letters. Since
$\I_0$ is a prime ideal, either $f\in \I_0$ or $g\in \I_0$.

Assume that $f\in \I_0$. Then there exists $F'\in \I_h$ with
$\rho(F')=f$ and $\deg(F')=\deg(f)\leq \deg(F)$. Denote $
F-F'=h\Delta F$; it is clear that $\deg(\Delta F)<\deg(F)$. If
$\Delta F\in\I_h$ then $F$ itself is in $\I_h$ and we are through;
otherwise observe that $$h\Delta F G =FG-F'G\in\I_h .$$ Since
$U_h/\I_h$ is without torsion (\cite{fl}), we can ``divide" by
$h$, and  it follows that $\Delta F G\in\I_h$.

We can now proceed to show that either $\Delta F$ or $G$ is in
$\I_h$. But notice that we have reduced the total degree. We can
apply the argument again until we arrive to the situation that one
of the factors has degree zero (it is a number). Then it follows
that the other factor is in $\I_h$ and eventually that $F$ or $G$
are in $\I_h$, as we wanted to prove. \hfill $\blacksquare$

\bigskip

We now want  to show that there exists two different
$\C[h]$-deformations on the same orbit that are not isomorphic. We
consider $\g=\rm{sl}_2(\C)$.  Let $$ I_h=(P-\mu^0), \qquad
I_h'=(P-\mu^0-\sqrt{2}h), $$ where $P$ is the quadratic Casimir.
Assume that $U_h/I_h \cong U_h/I_h'$. Since any isomorphism will
send the ideal $(h-1)$ into the ideal $(h-1)$, the quotient of
both algebras by $(h-1)$ must be isomorphic. But the   algebra
$U/(P-\mu^0)$ has   finite dimensional representations only for
certain values of $\mu^0$. In particular, for $\mu^0$ irrational,
it has no finite dimensional representations \cite{va}. It is
enough to take $\mu^0$ such that $U_h/I_h $ has finite dimensional
representations and we reach a contradiction.

The same is true for   formal deformations. In fact, with the same
reasoning as in Ref. \cite{va} we have that $U_{[h]}/(P-\mu(h))$
admits finite dimensional representations
 only for appropriate $\mu(h)=\mu^0$.

\section{Star products and equivalence}

\begin{definition}
Given  $\A_{[h]}$, a formal deformation of a Poisson algebra $\A$
and a $\C[[h]]$-module isomorphism
$\psi:\A[[h]]\longrightarrow \A_{[h]}$
as in Definition \ref{fd}, we say that the associative
product in $\A[[h]]$ defined by $$a\star
b=\psi^{-1}(\psi(a)\cdot\psi(b)), \qquad a,b\in \A[[h]] $$ is a
star product on $\A[[h]]$.
\end{definition}

It follows from property {\bf a} in Definition \ref{fd} that a
star product can always be written as
\begin{equation}
a\star b=ab +\sum_{n> 0}h^nB_n(a,b) \label{ps}
\end{equation}
where $B_n$ are bilinear operators and by juxtaposition $ab$ we
denote the commutative product in $\A[[h]]$.
Property {\bf b} in Definition \ref{fd} implies that
$\{a,b\}=B_1^-(a,b):=B_1(a,b)-B_1(b,a)$.

For a given $\A_{[h]}$ there are many choices of the isomorphism
$\psi$ (it is not canonical). Once $\psi$ is given, the star
product $\star$ is defined and one regards $\A[[h]]$ as an
associative non commutative $\C[[h]]$-algebra. Let $\star$ and
$\star'$ be different star products corresponding to the same
deformation, defined by the maps
\begin{xalignat*}{3}
&\psi:\A[[h]]\longrightarrow \A_{[h]},& &a\star
b=\psi^{-1}(\psi(a)\cdot\psi(b)),\\ &\psi':\A[[h]]\longrightarrow
\A_{[h]},& & a\star' b={\psi'}^{-1}(\psi'(a)\cdot\psi'(b)).
\end{xalignat*}
They define  isomorphic algebras. The isomorphism $T:\A[[h]]
\longrightarrow \A[[h]]$
 is  given by
$$ T={\psi'}^{-1}\circ \psi, \qquad T(a\star b)=T(a)\star' T(b).
$$ $T$ can also be expressed as a power series
\begin{equation}
T=\sum_{n\geq 0}h^nT_n. \label{iso}
\end{equation}
in terms of the linear operators $T_n$. It is easy to show that
$T_0$ is an automorphism of the commutative algebra $\A[[h]]$ $$
T_0(ab)=T_0(a)T_0(b),\ \ \ \ \ a,b\in\A[[h]] $$ and of the Poisson
algebra $\A$ $$ T_0\{a,b\}=\{T_0(a),T_0(b)\},\ \ \ \ \ a,b\in\A.
$$

\begin{definition}
If $\star$ and $\star'$ are two isomorphic star products on
$\A[[h]]$, the isomorphism being $T:\A[[h]] \longrightarrow
\A[[h]]$ as in (\ref{iso}), we say that they are gauge equivalent
if $T_0=\mbox{Id}$.
\end{definition}

A star product is {\it differential} if $\A=C^\infty(M)$ for a
smooth manifold  $M$,
 and  the operators $B_n$ in (\ref{ps}) are
bidifferential operators.  An example of differential star product
is the one induced on $\g^*$ by the map (\ref{sym}). It is in
principle defined on polynomials, but it can be extended to
$C^\infty(M)$ through operators $B_i$ that are bidifferential. It
was shown in Ref. \cite{de} that with a gauge transformation any
differential star product can be brought to a form under which the
bilinear operators $B_n$ are null on the constants (that is, the
zero degree doesn't appear).

 One can consider gauge
equivalence inside the class of differential star products by
considering only differential maps $T$. For this case, it was
shown by Kontsevich in \cite{ko} the following important theorem,

\begin{theorem}
\label{kontsevich} The set of gauge equivalence classes of
differential star products on a smooth manifold $M$ can be
naturally identified with the set of equivalence classes of
Poisson structures depending formally on $h$,
 $$ \alpha=h\alpha_1+h^2\alpha_2+\cdots $$
 modulo the action
of the group of formal paths in the diffeomorphism group of $M$,
starting at the identity isomorphism. \label{kt}
\end{theorem}
In particular, for a given Poisson structure $\alpha_1$,  we have
the equivalence class of differential star products canonically
associated to $h\alpha_1$.

\bigskip

 We
 explain briefly the concept of formal paths in the
diffeomorphism group of $M$. For further details we refer to Ref.
\cite{ko, ar}. Let $\m$ be the maximal ideal in $\R[[t]]$.
Consider $\La$ the algebra of  polyvector fields with the
Schouten-Nijenhuis bracket. It is a differential graded Lie
algebra with zero differential. We recall that a Poisson structure
is a bivector field such that its Schouten-Nijenhuis bracket with
itself is zero. Let $\La_0$ be the algebra of vector fields on
$M$. They are the  0-cochains of the  complex. Consider
$\La_0\otimes \m$. The exponential of this algebra is the group of
formal paths in the diffeomorphism group starting with the
identity. $\La_1$ is the set of (skew-symmetric) bivector fields.
$\La_0$   acts on $\La_1$ with the Schouten-Nijenhuis bracket,
\begin{eqnarray*} &Z(B)(f_1,f_2)= Z(B(f_1,f_2))
-B(Z(f_1),f_2)-B(f_1, Z(f_2)), \\&Z\in\La_0, \;\;B\in \La_1,\;\;
f_1,f_2\in C^\infty(M)
\end{eqnarray*}
 and this action can be exponentiated to the group.

\subsection{Uniqueness of the deformation of the Kirillov Poisson structure \label{uni}}

We want to determine whether the equivalence class of the Kirillov
Poisson bracket in $\g^*$ is the only class of formal Poisson
structures whose first order term is the Kirillov Poisson bracket
(as it happens, for example, in any flat symplectic manifold
\cite{bffls}). This is actually the case, at least for algebraic
Poisson structures (we say that a Poisson structure $\beta$ is
algebraic if $\beta(p,q)$ is a polynomial whenever $p$ and $q$ are
polynomials) and $\g$ semisimple.
\begin{proposition}Let $\g$ be a real semisimple algebra.
Let $\beta$ be an algebraic differential formal Poisson structure
$$\beta=\sum_{i=0}^{\infty}h^{i+1}\beta_i, $$
 such that  $\beta_0$ is the
Kirillov Poisson structure in $\g^*$.
Then $\beta$ is equivalent to $\beta_0$.
\end{proposition}
\noindent {\it Proof.}
The Jacobi identity at first
order is satisfied since $\beta_0$ itself is a Poisson structure and
at second order it implies that $\beta_1$ is a two-cocycle in the
Chevalley cohomology of $\beta_0$. If $\beta_1$ is a coboundary,
then $$\beta_1(f_1, f_2)=\delta Z(f_1,f_2)=Z(\beta_0(f_1,f_2))
-\beta_0(Z(f_1),f_2)-\beta_0(f_1, Z(f_2)), $$ with $Z$ a
0-cochain.
Then a gauge  transformation   (formal path in the
diffeomorphism group) of the form
 $$ \varphi=\mbox{Id} +h Z +\cdots $$
shows that $\beta$ is equivalent to a formal Poisson structure
without term of order $h^2$, {\em i.e.}, we can assume that
$\beta_1=0$. But then $\beta_2$ is a cocycle and we can proceed
recursively. Hence, to prove that $\beta\sim\beta_0$ it is
actually  sufficient to show that the Chevalley cohomology of
$\beta_0$ is zero. Since $\beta$ is a bivector field and it is
algebraic, it is sufficient to check that there is no  non trivial
algebraic two cocycle with order of differentiability (1,1).  We
will show that this is the case.

Such a cocycle is an antisymmetric bidifferential map, null on the
constants and with polynomial coefficients, $$
C^2:\mbox{Sym}(\g)\otimes\mbox{Sym}(\g)\longrightarrow M $$ where
$M=\mbox{Sym}(\g)$ is a left (Lie algebra) Sym$(\g)$-module, with
the action given by the Poisson bracket $\beta_0$. If $C^2$ has
order of differentiability (1,1), we can restrict $C^2$ non
trivially to first order polynomials. We denote that restriction
by
 $$\hat{C}^2:\g\otimes\g\longrightarrow\mbox{Sym}(\g). $$ Then $\hat
C^2$ is a cocycle in the  Lie algebra cohomology of order two of
$\g$ with values in Sym$(\g)$. Since $\g$ is semisimple, as a
consequence of Whitehead's lemma, this cohomology is zero (see for
example Ref. \cite{ka}). Hence $\hat C^2$ is trivial, {\em i.e.},
there exists a 1-cochain $\hat{C}^1:\g\longrightarrow
\mbox{Sym}(\g)$ such that $\hat{C}^2=\delta\hat{C}^1$. If
$\hat{C}^1$ is given on a basis of $\g$ by
$\hat{C}^1(X_i)=\hat{C}^1_i$, this means that $$ \hat
C^2(X_i,X_j)=\delta\hat{C}^1=\hat C_1([X_i,X_j])-\beta_0(X_i,\hat
C^1_j)-\beta_0(\hat C^1_i,X_j),$$ and $\hat C^1$ can be extended
to a 1-cochain in the Chevalley complex by $$C^1(f)=\hat
C^1_k\frac{\partial f}{\partial x^k}, \qquad f\in\mbox{Sym}(\g).$$
We then have $C^2=\delta C^1,$ showing that $C^2$ is trivial.
Hence $\beta$ is equivalent to $\beta_0$, as we wanted to show.
\hfill{$\blacksquare$}

\bigskip

Using Theorem \ref{kt}, we conclude that there is only one
equivalence class of star products whose first order term is the
Kirillov Poisson bracket.  All these star products  give algebra
structures on the polynomials on $\g^*$ isomorphic to $U_{[h]}$.

\section{Star products on the orbit}

A star product on the orbit  $\Theta$ associated to the
deformation of Theorem \ref{flt} is given by a linear isomorphism
$$ \tilde\psi:\C[\Theta][[h]]\longrightarrow U_{[h]}/\I_{[h]}. $$
In particular, if $\{x_{i_1}\cdots x_{i_k}, (i_1,\dots i_k)\in
S\}$ is a basis of $\C[\Theta]$ for some set of multiindices $S$,
then $\{X_{i_1}\cdots X_{i_k}, (i_1,\dots i_k)\in S\}$ is a basis
of $U_{[h]}/\I_{[h]}$ \cite{fl}. This defines a particular
isomorphism $\tilde\psi(x_{i_1}\cdots x_{i_k})=X_{i_1}\cdots
X_{i_k}$ and the corresponding star product.

This star product can be seen as the restriction to the orbit of a
star product on $\C[\g^*]$. We have only to extend the map
$\tilde\psi$ to an isomorphism $\psi: \C[\g^*][[h]]\longrightarrow
U_{[h]}$. This is guaranteed since $\C[\g^*]=\C[\Theta]\oplus
\I_0$, $U_{[h]}=U_{[h]}/\I_{[h]}\oplus \I_{[h]}$, and  $\I_0$ and
$\I_{[h]}$ are also isomorphic as $\C[[h]]$-modules.

We have then that the following diagram
\begin{equation}
\begin{CD}
\C[\g^*][[h]]@>\psi>>U_{[h]}\\ @VV{\pi}V @VV{\pi_h}V\\
\C[\Theta][[h]]@>\tilde\psi>>U_{[h]}/I_{[h]}
\end{CD}
\label{cd}
\end{equation}
commutes. In general, we say that a star product on $\g^*$ is {\it
tangential} to the orbit $\Theta$ if it defines a star product on
$\Theta$ by restriction. So the star product in (\ref{cd}) is
tangential.

\begin{example} Star product on an orbit of  {\rm SU(2)}.
\end{example}
 Consider the  Lie algebra of SU(2),
$$ [X,Y]=Z, \qquad [Y,Z]=X,\qquad [Z,X]=Y. $$ The subalgebra of
invariant polynomials on $\g^*$ is generated by $p=x^2+y^2+z^2$,
so the corresponding Casimir is $P=X^2+Y^2+Z^2$. We consider the
orbit $p=c^2$, $c\in \R,\, c\neq 0$. A basis of $\I_0$ is
$B_1=\{x^ry^sz^t(p-c^2),\; r,s,t=0,1,2,\dots\}$ and one can
complete it to a basis in $\g^*$ by adding $B_2=\{x^ry^sz^\nu,\;
\nu=0,1,\; r,s=0,1,2,\dots\}$. The equivalence classes of the
elements in $B_2$ are a basis of $\C[\g^*]/\I_0$.

Let $\I_{[h]}$ be the ideal in $U_{[h]}$  generated by $P-c^2$. We
can define the isomorphism $\psi:\C[\g^*][[h]]\longrightarrow
U_{[h]}$ as
\begin{eqnarray}
&\psi(x^ry^sz^t(p-c^2))=X^rY^sZ^t(P-c^2),\qquad
r,s,t=0,1,2,\dots\nonumber\\ &\psi(x^ry^sz^\nu)=X^rY^sZ^\nu,\qquad
\nu=0,1,\; r,s=0,1,2,\dots\label{psi}
\end{eqnarray}
Clearly $\psi(\I_0)=\I_{[h]}$, so the star product defined by $\psi$
is tangential to the orbit. It is easy to check that if we move to
a neighboring orbit, $p={c'}^2$, then $\psi$, as defined in
(\ref{psi}) doesn't preserve the new ideal, that is,
$\psi(\I'_0)\neq \I'_{[h]}$.
\bigskip

One can construct
 a star product that is tangential to all the orbits in a
neighborhood of the regular orbit (this is in fact the definition
of ``tangential star product" given in \cite{cgr}). If
$p_i=c^0_i$, $i=1,\dots m$ define the regular orbit
$\Theta_{(c^0_1,\dots c^0_m)}$ with ideal $\I_{(c^0_1,\dots
c^0_m)}$ one can construct a map $\psi$ such that
\begin{equation}
\psi(\I_{(c_1,\dots c_m)})=\I_{(c_1,\dots c_m),h} \label{prese}
\end{equation}
 for
$(c_1,\dots c_m)$ in a neighborhood of $(c^0_1,\dots c^0_m)$ and
$\I_{(c_1,\dots c_m),h}$ an ideal in $U_{[h]}$ of the type
required in Theorem \ref{flt}. The construction follows similar
lines to the one in  \cite{cg}. We consider the decomposition
$\C[\g^*]=I\otimes H$ where $I$ is the subalgebra of invariant
polynomials as in (\ref{invpol}) and $H$ is the set of harmonic
polynomials (this result is due to Kostant \cite{kos}). Harmonic
polynomials are  in one to one correspondence with the polynomials
on the orbit, so we have in fact $$\C[\g^*]\approx
I\otimes\C[\Theta_{(c^0_1,\dots c^0_m)}]. $$ Consider now the
basis in $I$ $\{(p_{i_1}-c_{i_1}) \cdots (p_{i_k}-c_{i_k}),
i_1\leq \cdots \leq i_k\}$ and the basis in
$\C[\Theta_{(c^0_1,\dots c^0_n)}]$ as before, $\{x_{j_1}\cdots
x_{j_l}, (j_1,\dots  j_l)\in S\}$. We define the $\C[[h]]$-module
isomorphism
\begin{equation}\label{tsp}
\begin{array}{l}
\psi\bigl((p_{i_1}-c^0_{i_1})\cdots (p_{i_k}-c^0_{i_k})\otimes
x_{j_1}\cdots x_{j_l}\bigr)=\\ \ \ \ \ \ \ \ \ \ \ \ \ \ \ \ \
=(P_{i_1}-c_{i_1}(h))\cdots (P_{i_k}-c_{i_k}(h))\otimes (X_{j_1}
\cdots X_{j_l}).
\end{array}
\end{equation}
 It is obvious that it preserves the ideal, $\psi(\I_{(c^0_1,\dots
c^0_n)})=\I_{(c^0_1,\dots c^0_m),h}$. A closer look reveals that,
in fact $\psi(\I_{(c_1,\dots c_m)})=\I_{(c_1,\dots c_m),h}$, and
then $\tilde \psi$ in (\ref{cd}) is well defined for any
$(c_1,\dots c_m)$ in a neighborhood of $(c^0_1,\dots c^0_m)$.
 Consequently we have a star product that is tangential to all the
orbits in a neighborhood of the regular orbit.

In \cite{cg} it is shown that for SU(2) a star product of this
type (with $c_{i_l}(h)=c_{i_l}^0$)  is not differential. More
generally, it was shown in \cite{cgr} the following theorem
\begin{theorem}
\label{ntsp} If $\g$ is a semisimple Lie algebra there is no
differential star product on any neighborhood of the origin in
$\g^*$ which is tangential to the coadjoint orbits. \label{tcgr}
\end{theorem}

The only property of tangential star products that is used in the
proof of this theorem is that if $f$ is a function that is
constant on the orbits (in particular, the quadratic Casimir
$p_1$), then, $g\star f=gf.$ It is easy to show that the
tangential star products defined by (\ref{tsp}) satisfy this
property on all $\g^*$ (in particular in a neighborhood of 0), so
they are not differential.

\bigskip
On  any regular Poisson manifold there exists a star product that
is tangential and differential \cite{ma}. But on all of $\g^*$,
which is not regular,  Theorem \ref{tcgr} states that a star
product with both properties does not exist. To induce a star
product on a particular orbit, it is enough to assume that the
star product on $\g^*$ is tangent to only such orbit. One can find
star products on $\rm{Pol}(\g^*)$ isomorphic to $U_h$ that
restrict well to only one orbit (in the sense of (\ref{cd})).
Example \ref{nd} shows one of such star products for $\g$=su(2).
We prove that it is not differential, so at least in this case,
the relaxation of the tangentiality condition does not allow in
general for differentiability. In Section 5 we will investigate
how these deformations are related to differential deformations on
the orbit.

\bigskip

\begin{example} Non differential star product on  $\g^*=\rm{su(2)}^*$\label{nd}.
\end{example}Consider again the Lie algebra of SU(2), with the same
notation, and the orbit $\Theta^r$ given by $p=1$. It is a
2-sphere in $\R^3$. Fix the star product $\star$ on $\Theta$ by
choosing the $\C[[h]]$-isomorphism
\begin{eqnarray*}
\tilde\psi:\C[\Theta][[h]] & \longrightarrow & U_{[h]}/\I_{[h]}
\\ x^ny^mz^\nu & \mapsto & X^nY^mZ^\nu,\qquad\nu=0,1,\;
m,n=0,1,2,\dots
\end{eqnarray*}
We regard the cartesian coordinates $x$ and $y$ as functions on
the sphere and let $V$ be an open set in $\Theta$ where $(x,y)$
are coordinates. On this open set $V$ the 1-forms $dx$ and $dy$
form a basis for the module of 1-forms. Let  $\dex$ and $\dey$ be
the elements of the dual basis, that is,   $\dex$ and $\dey$ are
vector fields on $V$ such that
$$\langle\dex,dx\rangle=\langle\dey,dy\rangle=1,
\qquad\langle\dey,dx\rangle=\langle\dex,dy\rangle=0.$$ Any
differential operator on $V$ is an element of the algebra
generated by functions and by $\dex$ and $\dey$. The advantage of
$\dex$ and $\dey$ is that they behave well on polynomials in $x$
and $y$. We have
 $$0=\dex(1)=\dex(x^2+y^2+z^2)=2x+2z\dex(z),$$
hence $\dex(z)=-\frac{x}{z}$ and $\dey(z)=-\frac{y}{z}$. Observe
that $\dex$ and $\dey$ commute.

Assume that $\star$ is differential, $$f\star g=\sum_{i\geq 0}
h^iB_i(f,g)$$ where $B_i$ are bidifferential operators. To
determine $B_i$ it is enough to compute them on the monomials $x$
and $y$. With the following lemma we compute $B_1$.

\begin{lemma} Let $p_1$, $p_2$ be two polynomials in
$x$ and $y$, then we have: $$ p_1 \star p_2 = p_1p_2-hz
\partial_y(p_1)\partial_x(p_2)\quad{\rm mod}(h^2) $$
\end{lemma}
\noindent {\it Proof.} It is enough to show it for $p_1$, $p_2$
monomials. Let $p_1=x^ny^m$, $p_2=x^ry^s$. We use induction on
$N=m+r$. For $N=0$ it is clear. Let $N>0$. By the definition of
$\star$, $$ p_1 \star p_2 =x^ny^m \star x^ry^s=x^n(y^m \star
x^r)y^s= x^n[(y^m \star x^{r-1}) \star x]y^s. $$ By induction we
have: $$ p_1 \star p_2 =x^n[(x^{r-1}y^m-hzm(r-1)y^{m-1}x^{r-2})
\star x]y^s\quad{\rm mod}(h^2), $$ and by induction again we have:
\begin{eqnarray*}
p_1 \star p_2 & = & x^n[x^{r}y^m-hzmy^{m-1}x^{r-1}
               -hzm(r-1)y^{m-1}x^{r-1}]y^s\quad{\rm mod}(h^2) \\
            & = & x^n[x^{r}y^m-hzmry^{m-1}x^{r-1}]y^s\quad{\rm mod}(h^2) \\
           & = &  x^{n+r}y^{m+s}-hzmry^{m-1+s}x^{r-1+n}\quad{\rm mod}(h^2),
\end{eqnarray*}
which is what we wanted to prove. \hfill $\blacksquare$

\bigskip

According to the previous lemma $$z\star z=
z^2-h\frac{xy}{z}\quad{\rm
mod}(h^2)=1-x^2-y^2-h\frac{xy}{z}\quad{\rm mod}(h^2) ,$$ on the
other hand, by definition, $$ z\star z=
\tilde\psi^{-1}(\tilde\psi(z)\tilde\psi(z))
=\tilde\psi^{-1}(Z^2)=1-x^2-y^2, $$ a contradiction that shows
that $\star$ cannot be differential.

\section{Algebraic and differential star products on the regular orbit.}

 Let us consider the regular orbit $\Theta^r$ as a symplectic manifold. By Theorem \ref{kt}
 we can associate to the Poisson structure a (equivalence class of) differential
star product. It was already known \cite{dl, fe} that a star
product exists for any symplectic  manifold. In fact, differential
star products are not in general unique. The space of equivalence
classes of differential star products  such that $$f\star g-g\star
f=h\{f,g\},$$ being $\{\,\cdot\,,\,\cdot\,\}$ a symplectic Poisson
bracket, are classified by the sequences $\{\omega_n\}_{n\geq 1}$
of de Rham cohomology classes in $H^2(M)$ such that $\omega_1$ is
the symplectic form associated to the Poisson bracket. In fact,
the space of equivalence classes of star products is a principal
homogeneous space under the group $H^2(M)[[h]]$ \cite{dl, de}.

The symplectic two form is not defined in arbitrary Poisson
manifolds, so the natural structure to consider is the Poisson
bivector. We want to describe the space of equivalence classes of
star products for symplectic manifolds in terms of the Poisson
bivector, being this approach closer to the one of Kontsevich's
theorem for arbitrary manifolds. Let $M$ be a symplectic manifold
and  consider $\omega_h=\sum_{j\geq 0}h^j\omega_j\in H^2(M)[[h]]$,
where $\omega_0$ is the original symplectic two form and
$\omega_j$ are closed two forms. Since $\omega_0$ is non
degenerate, $\omega_h$ defines an invertible map between tangent
and cotangent vector fields in the usual way,
$$\mu_h:\Gamma(TM)[[h]]\longrightarrow \Gamma(T^*M)[[h]],$$ which
can be extended to tensors. In fact, by closedness of $\omega_h$,
the map $$(f,g)\mapsto
\{f,g\}_h=h\omega_h(\mu_h^{-1}(df),\mu_h^{-1}(dg))$$ is a formal
Poisson structure in the sense of Kontsevich and this formal
Poisson structure is gauge equivalent to zero (the gauge group is
the group of formal paths in the diffeomorphism group starting
with the identity) if and only if all the $\omega_j$ are exact. We
have then that the set of equivalence classes in $H^2(M)[[h]]$ is
in one to one correspondence with the set of formal Poisson
structures modulo the action of the gauge group. So for symplectic
manifolds both descriptions, in terms of the symplectic form or in
terms of the Poisson bivector are equivalent.

Coadjoint orbits of compact groups are an example of manifolds
that admit inequivalent quantizations. In fact they have a
nontrivial de Rham cohomology $H^2(\Theta^r)$. In particular, the
symplectic form is a closed, non exact form, so we have many
inequivalent deformations.

 Let $\Theta^r$ be the   orbit defined by $$ p_i(x_1,\dots
x_n)=c^0_i,\quad c^0_i\in \R,\quad i=1,\dots m,$$ $\Theta^r$ is
regular if and only if the differentials $dp_i$ are independent.
One can consider all the regular orbits given by the constraints
$p_i=c_i$, $i=1,\dots m$ with  $(c_1,\dots c_m)$ in a neighborhood
of $(c^0_1,\dots c^0_m)$ where the differentials are still
independent. The set of these points is a neighborhood  $\N\approx
\Theta^r\times \R^m$ of the regular orbit. $\N$ is a regular
Poisson manifold. The Poisson structure on $\N$ can be seen as a
symplectic structure on $\Theta^r$ which depends on certain
parameters, the invariant polynomials  $p_i$, which determine the
leaf of the foliation in $\N$.

We now want to examine various star products  that can be defined
on the open set  $\N$. We can consider the star product $\star_S$
induced by $U_h$ by means of the map (\ref{sym}).  It is
differential, but not tangential. It was shown in Ref. \cite{ko}
that the canonical deformation of the Kirillov Poisson structure
on $\g^*$ is isomorphic to $U_h$.

We can also consider the  quantization of the Kirillov symplectic
structure on the orbit given by Kontsevich's theorem.    From the
local expression of Kontsevich's quantization, one can see that it
is a smooth with respect to the parameters $p_i$. Interpreting the
parameters as transverse coordinates, Kontsevich's theorem applied
on $\Theta^r$ gives indeed a star product on $\N$ that is
tangential and differential. We denote it by $\star_T$.

 Finally we can consider the star product $\star_P$ on $\N$,
induced by a map $\psi$ as in formula (\ref{tsp}). $\star_P$ is
tangential to the orbit, but, in general, not differential. To sum
up we get Table \ref{3sp}.

\begin{table}[ht]
\begin{center}
\begin{tabular} {|m{1cm}|m{5cm}|m{5cm}|}
\hline
 $*_S$  &Isomorphic to $U_h$ (on the polynomials), induced by Sym.&
 Defined on all $\g^*$, differential, not tangential.\\\hline
$*_P$ & Isomorphic to $U_h$, induced by a map $\psi$ like
(\ref{tsp}).& Defined on all $\g^*$, not differential (only given
on polynomials), tangential to the orbits.\\\hline $*_T$ &Gluing
Kontsevich construction on the leaves. &Defined on $\N$,
differential, tangential.
\\\hline
\end{tabular}
\caption{ Star products on $\N$}\label{3sp}
\end{center}
\end{table}

The relation among these star products on $\N$ and the
corresponding star products induced on the orbit $\Theta^r$ will
be studied in \cite{fl2}.
\bigskip

\section{Summary}

In this paper we consider different methods of quantization for
regular orbits of compact semisimple Lie groups. From the
algebraic point of view, one can obtain non isomorphic
deformations of the same Poisson structure. These deformations can
be compared with geometric quantization since the formulation is
in terms of a certain prime ideal in the enveloping algebra. The
comparison with differential deformations becomes more difficult
since the polynomials are ``global" objects, very different from
the ``local"  $C^\infty$ functions, and in fact we see that the
star products obtained are not differential in general. At the end
we define three  star products on a regularly foliated
neighborhood of the orbit.
\section*{Acknowledgements}

 We  want to thank E.
Arbarello for useful discussions.

\end{document}